\documentclass[8pt,a4paper]{article}
\usepackage{subfigure}
\usepackage{CJK,CJKnumb}
\usepackage{epsf, graphicx}
\usepackage{latexsym,amsfonts,amsbsy,amssymb}
\usepackage{graphicx}
\usepackage{amsmath}
\numberwithin{equation}{section}
\usepackage[thmmarks]{ntheorem}
\usepackage{float}
\usepackage{subfigure}

{
   \theoremstyle{nonumberplain}
   \theoremheaderfont{\bfseries}
   \theorembodyfont{\normalfont}
   \theoremsymbol{\mbox{$\Box$}}
   \newtheorem{proof}{Proof.}
}
\usepackage{color}
\newtheorem{theorem}{Theorem}[section]

\title{Bifurcation Analysis in a Continuous-Time Information Model with Discrete and Distributed Delays}
\author{Jingli Ren\textsuperscript{1} \footnote{Corresponding
author. E-mail: renjl@zzu.edu.cn. Tel:13653827917.},~~Fangzhi Yu\textsuperscript{1}
\\ \normalsize \textsuperscript{1}~School of Mathematics and Statistics, Zhengzhou University,
Zhengzhou 450001, China }

\date{}
\date{}
\begin{document}
\maketitle
\begin{center}
\begin{minipage}{13cm}
\noindent
$\displaystyle
\textbf{Abstract}$   In this paper, we consider a continuous-time model with discrete and distributed delays  to describe how two pieces of information interact in online social networks. Sufficient conditions are carried out to illustrate the stability of each equilibrium. Taking time delay as a bifurcation parameter, the system undergoes a sequence of Hopf bifurcation when this parameter passes through a critical value. By methods of multiple scales, we prove that the direction of Hopf bifurcation is depending on the condition which is related to delay.

 \noindent
$\displaystyle
\textbf{Keywords}$\ \ distributed delay, discrete delay, methods of multiple scales, Hopf bifurcation, direction, stability
 \end{minipage}
\end{center}
  \section{Introduction}

  In this paper, we study an information model with discrete and distributed delays in online social networks where two pieces of information interacting with each other when spreading.

   The proposed model is given as follows
  \begin{equation}
  \begin{cases}
  \begin{split}
  &\frac{{\rm d} u(t)}{{\rm d} t}={r}_{1}{u}(t)({1}-{a}_{1}{u}(t))-b_{1}r_{1}{u}(t-s){v}(t-s), \\
  &\frac{{\rm d}{v}(t)}{{\rm d} t}=r_{2}{v}(t)({1}-{a}_{2}{v}(t))+\int_0^\infty{b}_{2}r_{2}{f}(\tau)e^{-\mu \tau}{u}(t-\tau){v}(t-\tau)d\tau.
  \end{split}
  \end{cases}
  \end{equation}

   In the model, the variable ${u}(t)$ and ${v}(t)$ denote the density of influenced users for two pieces of information at time $t$ respectively, they represent the ratio of the number of influenced users at time $t$ over the total number of users. ${r}_{1}$ and $r_2$ denote the intrinsic growth  rates which indicates how fast the information propagates and they are positive constants. ${1}/{a_1}>0$ and ${1}/{a_2}>0$ are the maximum carrying capacities. To incorporate delays into the model, we briefly summarize the main stage of the process. We will denote two pieces of information with $u$ and $v$ for convenience. Without ${v}$, ${u}$ admits the logistic equation $\dot{u}(t)={r}_{1}{u}(t)({1}-{a}_{1}{u}(t))$. We take ${s}$ as the time delay between ${v}$ contacts ${u}$ and take effect. The intervening rate of ${v}$ to ${u}$ is ${b}_{1}\in\mathbb{R}$. Once ${v}$ presents, there will be an effective term which is related to the number of two information at ${s}$ time ago. Similarly, in the absence of ${u}$, ${v}$ also follows a logistic equation. In the presence of ${u}$,we consider $\tau$ as the random variable that describes the time between information ${u}$ bump with ${v}$ and influence with a probability distribution $f(\tau)$ and ${f}(\tau):[{0},\infty)\rightarrow[0,\infty)$ is probability distribution with compact support. ${b}_{2}\in\mathbb{R}$ is the intervening rate of ${u}$ to ${v}$. The presence of the distributed time delays must not affect the equilibrium values, so we normalize the kernels such that $\int_0^\infty{f}(\tau)d{\tau}=1$, ${f}(\tau)\geq{0}$. The density of influenced ${v}$ at time $t$ is given by $$\int_0^\infty{b}_{2}{r}_{2}{f}(\tau){e}^{-\mu \tau}{u}(t-\tau){v}(t-\tau)d\tau,$$ where ${e}^{-\mu \tau}$ accounts for the loss rate of ${v}$ during time period $[t-\tau, t]$ due to ${u}$.

 Social networking software, for example Facebook and Wechat, have been extremely important ways for users to get information and communicate. Online social networks provide new insights into human social behavior  that researchers show great interests in this filed. If there is only one piece of information in the social network, the dissemination pattern is logistic as the following figure.
 \begin{figure}[H]
\centering
\label{Fig.sub.1}
\includegraphics[width=0.7\textwidth]{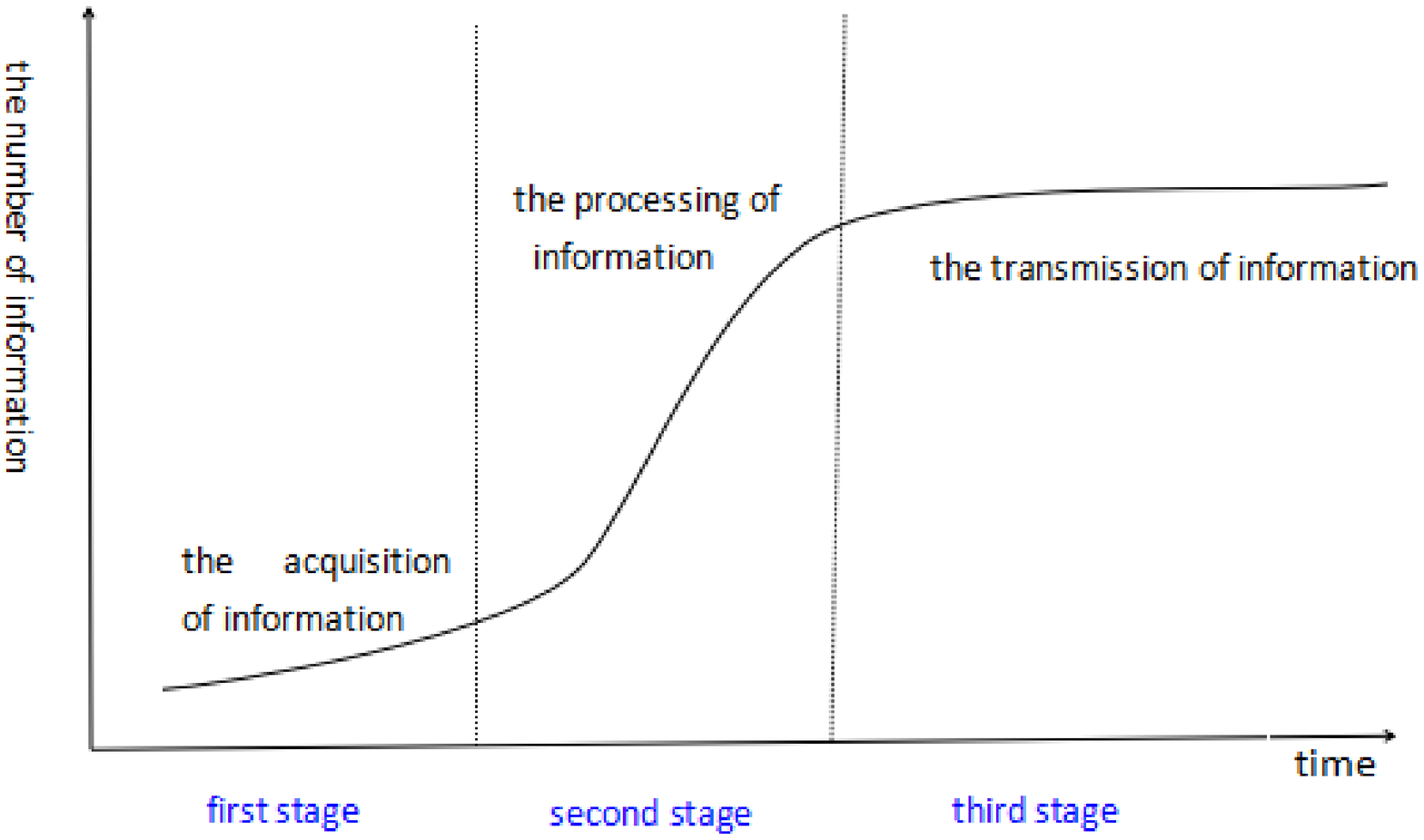}
\end{figure}

  Based on the mechanism and characteristic of information transmission, F. Wang, H.Y. Wang and K. Xu divided the process into two stages: growth process and social process in \cite{Hy}. They proposed a temporal and spatial pattern logistic PDE model to illustrate how one pieces of information diffuse in online social network. They also gave their prediction with the real data from Digg. with a high accuracy. In \cite{Lei}, C.X. Lei, Z.G. Lin and H.Y. Wang took the model as a free boundary problem to have a further research. They provided the condition where the information is spreading or vanishing. Considering the complexity of social network, C. Peng, K. Xu, F. Wang and H.Y. Wang started to consider information from multiple sources by analysing the real data from Digg in \cite{Cp}. In \cite{Freemana}, M. Freemana, J.M. Vittie, I. Sivak and J.H. Wu come up with using SIR model to describe the information diffusion. In \cite{Chunru}, C.R. Li and Z.J. Ma discussed the Hopf bifurcation of the above model with a discrete delay. In \cite{Zhu}, L.H. Zhu, H.Y. Zhao and H.Y. Wang proposed a SI model with a set of coupled PFDEs to describe spatial-temporal dynamics of rumor propagation. In \cite{Ren}, J.L. Ren and L.P. Yu constructed a codimension-two discrete model to study two pieces of information with interaction. They classified the model into three types: reverse type, intervention type and mutualistic type. Then they discussed 1:2, 1:3, 1:4 resonances respectively and give two control strategies. We will further discuss the continuous-version model in this paper.

  On theory, our model has not been researched before. There is a lot of papers about the classic Lokta-Volterra model, which is similar to our model. Since it is reasonable to incorporate time delays to Lokta-Volterra model, there are extensive research in this field.
  At first, researchers consider only one single time delay. In \cite{W}, W. Aiello and H. Freedman considered the time between immature stage and mature stage as a constant delay and construct a model of single species. They analysed the stability of equilibrium by computing the characteristic equation and proved the model has nonoscillatory solutions.  Song and Wei \cite{Song} considered a delayed predator-prey system to study the existence of local Hopf bifurcations and derived explicit formulas to determine the stability and the direction of periodic solutions bifurcating from Hopf bifurcations by using the normal form theory and center manifold argument. M. Bandyopadhyay and Sandip Banerjee discussed the time interval between the moments that an individual prey is killed and the corresponding biomass is added to the predator population as a parameter to construct a more compliticated stage-structured predator-prey model \cite{M}. They believed that the time delays have an significant effect to the model and will lead to oscillations in constant environment. For the sake of have a more generalized and accurate application to reality, researchers  begin to deliberate two time delays. In \cite{Faria}, Faria introduced two different delays (hunting delay and delay in the predator maturation) into the system and took one of the delays as a parameter to study the Hopf bifurcation when it crosses some critical values
\begin{equation*}
\begin {cases}
\dot{x}(t)=x(t)[r_1-a_{11}x(t)-a_{12}y(t-\tau_1)],\\
\dot{y}(t)=y(t)[-r_2+a_{21}x(t-\tau_2)-a_{22}y(t)].
\end {cases}
\end{equation*}
Xu \emph{et.al} considered a more complex situation and investigated the following delayed Lokta-Volterra model with two delays just as Faria mentioned above
\begin{equation*}
\begin {cases}
\dot{x}(t)=x(t)[r_1-a_{11}x(t-\tau_1)-a_{12}y(t-\tau_2)],\\
\dot{y}(t)=y(t)[-r_2+a_{21}x(t-\tau_1)-a_{22}y(t-\tau_2)].
\end {cases}
\end{equation*}
 The linear stability of the above system was investigated, Hopf bifurcation and the direction of bifurcating periodic solutions were demonstrated in \cite{Xu}.
In \cite{Yuan}, Song and Yuan considered the following model with a discrete time delay and a distributed time delay of the form
\begin{equation*}
\begin {cases}
\dot{x}(t)=x(t)[r_1-a_{11}x(t)-a_{12}y(t-\tau)],\\
\dot{y}(t)=y(t)[-r_2+a_{21}\int_{-\infty}^tF(t-s)x(s)ds-a_{22}y(t)].
\end {cases}
\end{equation*}
They considered the stability of the positive equilibrium and the existence of local Hopf bifurcations and derived the explicit formulas determining stability, direction and other properties of bifurcating periodic solutions.

  In this paper, we let $f(\tau)$ be the exponential kernels $f(\tau)=re^{-r\tau}(r>0)$, ${w}=\int_0^\infty{f}(\tau)e^{-\mu \tau}{u}(t-\tau){v}(t-\tau)d\tau$, then system (2.1) can be transformed into
\begin{equation}\label{equation:eq-a1}
\begin{cases}
\begin{split}
&\frac{{\rm d} {u}(t)}{{\rm d}t}={r}_{1}{u}(t)({1}-{a}_{1}{u}(t))-b_{1}r_{1}{u}(t-s){v}(t-s),\\
&\frac{{\rm d} {v}(t)}{{\rm d}t}=r_{2}{v}(t)({1}-{a}_{2}{v}(t))+{b}_{2}r_{2}{w}(t),\\
&\frac{{\rm d} {w}(t)}{{\rm d}t}={u}(t){v}(t)-(\mu +r){w}(t).
\end{split}
\end{cases}
 \end{equation}

 The rest of the paper is organized as follows.
In section 2, we compute the equilibria and investigate their stability respectively.
In section 3, we obtain the critical value that  Hopf bifurcation occur. In section 4, the specific form and conditions determining the direction of the local Hopf bifurcation are given. Simulations are carried out to illustrate the theory results in section 5.

\section{The stability of Equilibrium}
    From (1.2), we can easily obtain four equilibria ${E}_{0}(0,0,0)$, ${E}_{1}({\frac{1}{a_{1}}},{0},0)$,
  \noindent
$\displaystyle$
    ${E}_{2}({0},{\frac{1}{a_{2}},0})$ and ${E}_{\ast}({u}_{\ast},{v}_{\ast},{w}_{\ast})$ where

\begin{equation*}
  \begin{cases}
  \begin{split}
  &{u}_\ast=\frac{(a_{2}-b_{1})(\mu +r)}{a_{1}a_{2}(\mu +r)+b_{1}b_{2}},\\
  &{v}_\ast=\frac{a_{1}(\mu+r)+b_{2}}{a_{1}a_{2}(\mu+r)+b_{1}b_{2}},\\
  &{w}_\ast=\frac{(a_{2}-b_{1})[a_{1}(\mu+r)+b_{2}]}{[a_{1}a_{2}(\mu+r)+b_{1}b_{2}]^{2}}.\\
  \end{split}
  \end{cases}
\end{equation*}
Next we will explore the stability of four equilibrium.

 The characteristic equation of system (1.2) at ${E}_{0}(0,0,0)$ is $$\varphi(\lambda)=(\lambda-r_{1})(\lambda-r_{2})(\lambda+\mu+r).$$
  Obviously that the characteristic equation has at least two positive roots $r_{1}$ and $r_{2},$ which means  ${E}_{0}(0,0,0)$  is unstable.

Similarly, we can get the characteristic equation  at ${E}_{1}(\frac{1}{a_{1}},0,0)$ is $$\varphi(\lambda)=(\lambda+r_{1})[\lambda^{2}+(\mu+r-r_{2})\lambda-(\mu+r)r_{2}-\frac{b_{2}r_{2}}{a_{1}}].$$

 If  $\mu+r+\frac{b_{2}}{a_{1}}>{0}$, it is clear that
 $\varphi(0)<{0}$ and $\varphi(\lambda)_{\lambda\rightarrow+\infty}\rightarrow+\infty$. This yields that the equation above has at least one positive root, which implies that the equilibrium ${E}_{1}(\frac{1}{a_{1}},0,0)$ is unstable.
\begin{theorem}
   For system (1.2),\\
   (1)if $b_{1}>a_{2}$$(b_{1}<a_{2})$, ${E}_{2}(0,\frac{1}{a_{2}},0)$ is locally asymptotically stable (unstable).\\
   (2)if and only if $b_{1}<a_{2}$ and $a_1a_2(\mu+r)>\max\{-b_1b_2, -a_2b_2\}$ are satisfied at the same time, there will be a positive equilibrium ${E}_\ast({u}_\ast, {v}_\ast, {w}_\ast)$.
\end{theorem}

\begin{proof}
The characteristic equation of system (1.2) at ${E}_{2}(0,\frac{1}{a_{2}},0)$ is $$\varphi(\lambda)=(\lambda-r_{1}+\frac{b_{1}r_{1}}{a_{2}})(\lambda+r_{2})(\lambda+\mu+r).$$

 It is obvious that when $b_{1}>a_{2},$ the equation above has three negative roots and ${E}_{2}$ is locally asymptotically stable. When $b_{1}<a_{2}$ , the equation has a positive root, which means the system is unstable.
 From the form of ${E}_\ast({u}_\ast,{v}_\ast,{w}_\ast)$, we can know that positive root ${E}_\ast({u}_\ast,{v}_\ast,{w}_\ast)$ exists if and only if conditions $b_{1}<a_{2}$ and $a_1a_2(\mu+r)>\max\{-b_1b_2,-a_2b_2\}$ are satisfied.
 \end{proof}
 Remark: if $b_{1}>a_{2}$, ${E}_{2}(0,\frac{1}{a_{2}},0)$ is locally asymptotically stable, it means the density of users will not change again with the time at last.
\section{The Existence of Hopf Bifurcation}

  Let ${\bar{u}}={u}-{u}_\ast$, ${\bar{v}}={v}-{v}_\ast$, ${\bar{w}}={w}-{w}_\ast$, the system (1.2) can be transformed into the following form
\begin{equation}
\begin{cases}
\begin{split}
  &\frac{{\rm d} {\bar{u}}}{{\rm d}t}=r_{1}{\bar{u}}(1-2a_{1}{u}_\ast)-b_{1}r_{1}{v}_\ast{\bar{u}}(t-s)-b_{1}r_{1}{u}_\ast{\bar{v}}(t-s),\\
  &\frac{{\rm d}{\bar{v}}}{{\rm d}t}=r_{2}{\bar{v}}(1-2a_{2}{v}_\ast)+b_{2}r_{2}{\bar{w}}, \\
&\frac{{\rm d} {\bar{w}}}{{\rm d}t}={u}_\ast{\bar{v}}+{v}_\ast{\bar{u}}-(\mu+r){\bar{w}}.
\end{split}
\end{cases}
\end{equation}
Replace ${\bar{u}} ,{\bar{v}}, {\bar{w}}$ with ${u}, {v}$ and  ${w}$, the system can be written as
\begin{equation}
\begin{cases}
\begin{split}
  &\frac{{\rm d}{u}}{{\rm d}t}=r_{1}{u}(1-2a_{1}{u}_\ast)-b_{1}r_{1}{v}_\ast{u}(t-s)-b_{1}r_{1}{u}_\ast{v}(t-s),\\
  &\frac{{\rm d}{v}}{{\rm d}t}=r_{2}{v}(1-2a_{2}{v}_\ast)+b_{2}r_{2}{w}, \\
  &\frac{{\rm d}{w}}{{\rm d}t}={u}_\ast{v}+{v}_\ast{u}-(\mu+r){w}.
  \end{split}
  \end{cases}
\end{equation}
The associated characteristic equation at ${E}_\ast$ is
\begin{equation}\label{}
  \lambda^3+p_{2}\lambda^2+p_{1}\lambda+p_{0}+(q_{2}\lambda^2+q_{1}\lambda+q_{0})e^{-\lambda s}=0,
\end{equation}
where
\begin{eqnarray*}
  && p_{2}=\mu+r-r_{2}-r_{1}+2a_{1}r_{1}{u}_\ast+2a_{2}r_{2}{v}_\ast,\\
  && p_{1}=(2a_{1}r_{1}{u}_\ast-r_{1})(\mu+r-r_{2}+2a_{2}r_{2}{v}_\ast)-(\mu+r)(r_{2}-2a_{2}r_{2}{v}_\ast)-b_{2}r_{2}{u}_\ast, \\
  &&p_{0}=(r_{1}-2a_{1}r_{1}{u}_\ast)[(\mu+r)(r_{2}-2a_{2}r_{2}{v}_\ast)+b_{2}r_{2}{u}_\ast],  \\
  &&q_{2}=b_{1}r_{1}{v}_\ast,  \\
  &&q_{1}=b_{1}r_{1}{v}_\ast(\mu+r-r_{2}+2a_{2}r_{2}{v}_\ast),  \\
  &&q_{0}=-b_{1}r_{1}{v}_\ast(\mu+r)(r_{2}-2a_{2}r_{2}{v}_\ast).
\end{eqnarray*}

When $s={0},$ the characteristic equation becomes
 \begin{equation}
   \label{equation:eq-a1}
 \lambda^3+(p_{2}+q_{2})\lambda^2+(p_{1}+q_{1})\lambda+(p_{0}+q_{0})=0.
\end{equation}

\noindent
$\displaystyle$
Firstly, we make a hypothesis as follow

  $p_{0}+q_{0}>0$,\ \ \ \ $(p_{2}+q_{2})(p_{1}+q_{1})>q_{0}+p_{0}$\ \ \ \ \ \ \ \ \ \ \ \ \ \ \ \  \ (H1)

When (H1) is satisfied, all roots of eq. (3.4) have negative real parts and ${E}_\ast$ is locally asymptotically stable by using the Routh-Hurwitz criterion.

 Let $\lambda=i\omega(\omega>0)$ be the purely imaginary root of eq. (3.3). Then $\omega$  satisfies the following equation
 \begin{equation*}
   -i\omega^3-p_2\omega^2+p_1\omega i+p_0+(q_2\omega^2i+q_1\omega i+q_0)(cos\omega s-isin\omega s)=0.
 \end{equation*}

 Separating real and imaginary parts respectively, we have
 \begin{eqnarray}
\label{equation:eq-a6}
   &&q_{1}\omega\sin(\omega s)+(q_{0}-q_{2}\omega^2)\cos(\omega s)=p_{2}\omega^{2}-p_{0}, \\
   \label{equation:eq-a7}
   &&q_{1}\omega\cos(\omega s)-(q_{0}-q_{2}\omega^2)\sin(\omega s)=\omega^3-p_{1}\omega.
   \end{eqnarray}

  Squaring and adding the above equation, we will get
   \begin{equation}
     \omega^6+m\omega^4+n\omega^2+h=0,
   \end{equation}
 where
   \begin{eqnarray*}
     &&m=p_{2}^2-q_{2}^2-2p_{1}, \\
     &&n=p_{1}^2+2q_{0}q_{2}-q_{1}^2-2p_{0}p_{2},\\
     &&h=p_{0}^2-q_{0}^2.
   \end{eqnarray*}

 \noindent
$\displaystyle$
  Denote $z=\omega^2$, eq. (3.7) can be expressed by

   $$G(z)=z^3+mz^2+nz+h=0.$$

It is obvious that if $h<0$, then $G(0)<0$ and $\lim\limits_{z\rightarrow\infty}G(z)\rightarrow+\infty$. So  $G(z)$ has at least one positive real root.

   Without loss of generality, we assume that it has three positive roots: $z_{1}$, $z_{2}$ and $z_{3}$. Then $\omega_1=\sqrt{z_1}$, $\omega_2=\sqrt{z_2}$ and $\omega_3=\sqrt{z_3}$. From (3.5)-(3.6), we can get
\begin{equation*}
  \begin{split}
  s_{k}^{(j)}=\frac{1}{\omega_{k}}[arccos\frac{(q_{0}-q_{2}\omega_{k}^2)(p_{2}\omega_{k}^2-p_{0})+q_{1}\omega_{k}(\omega^3-p_{1}\omega_{k})}{(q_{0}-q_{2}\omega_{k}^2)^2+q_{1}^2\omega_{k}^2}+2j\pi], k=1,2,3;j \in N.
  \end{split}
\end{equation*}

   When $s=s_{k}^{(j)}$, $\pm i \omega$ are a pair of purely imaginary roots of the characteristic equation.

   Substituting $\lambda(s)$ into  eq. (3.3) and taking derivation with respect to $s$, we have
\begin{equation}
  (\frac{{\rm d}\lambda}{{\rm d}s})^{-1}=\frac{(3\lambda^2+2p_{2}\lambda+p_{1})e^{\lambda s}}{\lambda(q_{2}\lambda^{2}+q_{1}\lambda+q_{0})}+\frac{2q_{2}\lambda+q_{1}}{\lambda(q_{2}\lambda^2+q_{1}\lambda+q_{0})}-\frac{s}{\lambda},\\
  \end{equation}

   then
\noindent
$\displaystyle$
  \begin{equation*}
  \begin{split}
  Re(\frac{{\rm d}\lambda}{{\rm d}s})^{-1}_{s=s^j_k}&=Re\{\frac{(3\lambda^2+2p_{2}\lambda+p_{1})e^{\lambda s}}{\lambda(q_{2}\lambda^{2}+q_{1}\lambda+q_{0})}\}|_{s=s^j_k}+\frac{2q_{2}(q_{0}-q_{2}\omega^2)-q_{1}^{2}}{(q_{0}-q_{2}\omega^{2})^2+q_{1}^{2}\omega^{2}}\\
  \end{split}
  \end{equation*}
 \begin{equation*}
   \begin{split}
   &=Re{\frac{(-3\omega^2+2p_{2}\omega i+p_{1})(cos\omega s+isin\omega s)}{-q_{2}\omega^{3}i-q_{1}\omega^2+q_{0}\omega i}}|_{s=s^j_k}+\frac{2q_{2}(q_{0}-q_{2}\omega^2)-q_{1}^{2}}{(q_{0}-q_{2}\omega^{2})^2+q_{1}^{2}\omega^{2}}\\
   &=Re{\frac{(3\omega^2-2p_{2}\omega i-p_{1})(cos\omega s+isin\omega s)}{q_{1}\omega^{2}i-(q_{0}\omega-q_{2}\omega^3 )i}}|_{s=s^j_k}+\frac{2q_{2}(q_{0}-q_{2}\omega^2)-q_{1}^{2}}{(q_{0}-q_{2}\omega^{2})^2+q_{1}^{2}\omega^{2}}\\
   &=\frac{1}{[q_{1}^2\omega^2+(q_{0}-q_{2}\omega^2)]\omega^2}(3\omega^4+2m\omega^2+n)\\
   &=\frac{1}{[q_{1}^2\omega^2+(q_{0}-q_{2}\omega^2)]\omega^2}G'(z).
   \end{split}
 \end{equation*}

If $G'(z)\neq0$, the transversality condition that Hopf bifurcation can occur is satisfied.

Define $s_{0}=\min\limits_{k=1,2,3}\{s_{k}^{(0)}\}$, we can get the following theorem.

$\mathbf{Therom\  3.1}$ For system (1.2), suppose that $h<0$ and $G(z)\neq0$, then following results hold true.

(1)The system is asymptotically stable at $E_\ast$ for $s\in[0,s_{0})$.

(2)The system undergoes a Hopf bifurcation at the positive equilibrium when $s=s^j_k, k=1,2,3; j\in N$.

Remark: the time delay will cause period oscillation at the positive equilibrium point when it cross some critical value, which implies the information explosively spreads up and down in a short period.
  \section{Direction of Hopf bifurcation}

From the above discussion, we can get Hopf bifurcation at the critical $s_\ast=s^j_k$. Next in order to get the direction of  Hopf bifurcation, we discuss $s=s_\ast+\varepsilon^2\delta$ near $s_\ast$ and $\varepsilon$ is a nondimensional book-keeping parameter. In this paper, we apply multiple scales methods to solve this problem \cite{Peng}.

Moving the unique equilibrium to origin and rescaling $t \rightarrow s\bar{t}$, we can get
\begin{equation}
  \frac{{\rm d}X}{{\rm d}\bar{t}}=sAX+sA_sX_s+sF.
\end{equation}

Dropping the bars for simplification of notations, the system becomes
 \noindent
$\displaystyle$
\begin{equation}
  \frac{{\rm d}X}{{\rm d}t}=sAX+sA_sX_s+sF.
\end{equation}
where
$X=({u}(t),{v}(t),{w}(t))^T$ and $X_s=X(t-1)=({u}(t-1),{v}(t-1),{w}(t-1))^T,$

\begin{equation*}
A=\left(
 \begin{array}{ccc}
    r_1(1-2a_1u_\ast) &0 & 0 \\
    0 & r_2(1-2a_2v_\ast) & b_2r_2 \\
    v_\ast &u_\ast  & -(\mu+r) \\
  \end{array}
\right),\\
\end{equation*}

\begin{equation*}
A_s=\left(
 \begin{array}{ccc}
    -b_1r_1v_\ast &-b_1r_1u_\ast & 0 \\
    0 & 0 & 0\\
    0 &0 & 0 \\
  \end{array}
\right),
\end{equation*}

\begin{equation*}
F=\left(
 \begin{array}{ccc}
   -a_1r_1u^2-b_1r_1{u}(t-s){v}(t-s)\\
-a_2r_2{v}^2\\
{u}{v}
  \end{array}
\right).\\
\end{equation*}

We introduce the new time scale $T_0=t$ and $T_2=\varepsilon^2t$ . The solution does not depend on $T_1=\varepsilon t$ because secular terms first appear at $O(\varepsilon^3)$. We seek a uniform second-order approximate solution of system (4.2) in the form
\begin{equation}
  X(t,\varepsilon)=X(T_0,T_2)=\sum_{i=1}^{\infty}\varepsilon^iX_i(T_0,T_2).
\end{equation}

 \noindent
$\displaystyle$
The derivative with respect to $t$ is

\begin{equation}
  \frac{\rm d}{{\rm d}t}=\frac{\partial}{\partial T_0}+\varepsilon^2\frac{\partial}{\partial T_2},
\end{equation}

$X_s$ can be expressed as
\begin{equation}
  X(t-1,\varepsilon)=X(T_0-1,T_2-\varepsilon^2)=\sum_{i=1}^{3}\varepsilon^iX_{is}(T_0-1,T_2)-\varepsilon^3X_{1s}(T_0-1,T_2)+\cdot\cdot\cdot .
\end{equation}

Substituting eq. (4.3)-(4.5) into (4.2), comparing coefficients of  powers of $\varepsilon$, we can get
 \noindent
$\displaystyle$
\begin{eqnarray}
&D_0X_1-s_\ast AX_1-s_\ast A_sX_{1s} = 0, \\
&D_0X_2-s_\ast AX_2-s_\ast A_sX_{2s} = s_\ast P, \\
&D_0X_3-s_\ast AX_3-s_\ast A_sX_{3s}+D_2X_1-\delta(AX_1+A_sX_{1s})+s_\ast A_sD_2X_{1s} = s_\ast Q,
\end{eqnarray}

 \noindent
$\displaystyle$
where
\begin{equation*}
P=
\left(
\begin{matrix}
   -a_1r_1{u}_1^2-b_1r_1{u}_{1s}{v}_{1s}\\
-a_2r_2{v}_1^2\\
{u}_1{v}_1
\end{matrix}
\right ),
\end{equation*}
\begin{equation*}
Q=
\left(
\begin{matrix}
   -2a_1r_1{u}_1{u}_2-b_1r_1{u}_{1s}{v}_{2s}-b_1r_1{u}_{2s}{v}_{1s}\\
-2a_1r_1{v}_1{v}_2\\
{u}_1{v}_2+{u}_2{v}_1
\end{matrix}
\right ).
\end{equation*}

 \noindent
$\displaystyle$

 The solution of eq.(4.6) can be expressed as
\begin{equation*}
\begin{split}
  X_1=\mathbf{c}e^{i\omega_\ast s_\ast T_0}H(T_2)+\mathbf{\bar{c}}e^{-i\omega_\ast s_\ast T_0}\bar{H(T_2)},
  \end{split}
\end{equation*}
with $\mathbf{c}=(c,1,d)^{T}$ where
\begin{eqnarray*}
\begin{split}
&c=\frac{(i\omega_\ast+\mu+r)(i\omega_\ast-r_2+2a_2r_2{v}_\ast)-b_2r_2{v}_\ast}{{v}_\ast},  \\
&d=\frac{i\omega_\ast-r_2+2a_2r_2{v}_\ast}{b_2r_2}.
\end{split}
\end{eqnarray*}

 \noindent
$\displaystyle$
$X_2$ has a particular solution

$X_2=H^2(T_2)e^{2i\omega\ast s\ast T_0}\mathbf{e}+H(T_2)\bar{H(T_2)}\mathbf{f}+\bar{H(T_2)}^2e^{-2i\omega\ast s\ast T_0}\mathbf{\bar{e}}.$

Substituting it into eq. (4.7) and comparing the coefficients of\  $e^{2i\omega\ast s\ast T_0}H^2(T_2)$\ and \ $H(T_2)\bar{H(T_2)}$, we can obtain

\begin{equation*}
\mathbf{e}=
\left(
\begin{matrix}
   \xi-\sigma e_2\\
e_2\\
\frac{a_2}{b_2}+\eta e_2
\end{matrix}
\right ),\ \ \
\mathbf{f}=
\left(
\begin{matrix}
   m_1+n_1f_2\\
f_2\\
\frac{2a_2}{b_2}+pf_2
\end{matrix}
\right ),
\end{equation*}
  where
\begin{eqnarray*}
\begin{split}
   & \xi=\frac{-a_1r_1c^2-b_1r_1ce^{-2i\omega_\ast s_\ast}}{2i\omega_\ast-r_1+2a_1r_1{u}_\ast}+b_1r_1e^{-2i\omega_\ast s_\ast{v}_\ast}, \\
  & \sigma=\frac{b_1r_1e^{-2i\omega_\ast s_\ast}}{2i\omega_\ast-r_1+2a_1r_1{u}_\ast}+b_1r_1e^{-2i\omega_\ast s_\ast{v}_\ast}, \\
 & \eta=\frac{2i\omega_\ast-r_2+2a_2r_2{v}_\ast}{b_2r_2},\\
 &e_2=\frac{b_2\xi {v}_\ast+b_2c-a_2(2i\omega_\ast+\mu+r)}{2i\omega_\ast\eta+\mu \eta+r\eta+\sigma
  {v}_\ast-{u}_\ast},\\
   &m_1=\frac{2a_1c\bar{c}+b_1c+b_1\bar{c}}{1-2a_1{u}_\ast-b_1{v}_\ast},  \\
  &n_1=\frac{b_1{u}_\ast}{1-2a_1{u}_\ast-b_1{v}_\ast}, \\
  & p=\frac{2a_2{v}_\ast-1}{b_2},\\
  &f_2=\frac{2a_2(\mu+r)-b_2c-b_2\bar{c}-{v}_\ast m}{b_2{v}_\ast n+b_2{u}_\ast-\mu b_2p-rb_2p}.
  \end{split}
\end{eqnarray*}

 \noindent
$\displaystyle$
Substituting $X_1$ and $X_2$ into eq. (4.8), we obtain
\begin{equation*}
  D_0X_3-s_\ast AX_3-s_\ast A_sX_{3s}=-e^{i\omega_\ast s_\ast T_0}H'(\mathbf{I}+s_\ast A_se^{-i\omega_\ast s_\ast})\mathbf{c}+\delta e^{i\omega_\ast s_\ast T_0}H(A+A_se^{-i\omega_\ast})\mathbf{c}\\
\end{equation*}
\begin{equation}
  +s_\ast Me^{i\omega_\ast s_\ast T_0}H^2\bar{H}+cc+ST,
\end{equation}
with
\begin{equation*}
M=
\left(
\begin{matrix}
  -b_1r_1(cf_2+\bar{c}e_2+m_1+n_1f_2+\xi-\sigma e_2)e^{-i\omega_\ast s_\ast}\\
0\\
 cf_2+\bar{c}e_2+m_1+n_1f_2+\xi-\sigma e_2
\end{matrix}
\right ).
\end{equation*}

In the above equation,  $cc$ represents the complex conjugate terms and $ST$ denotes terms that do produce secular terms.

We can substitute $X_3=e^{i\omega_\ast s_\ast T_0}E(T_2)+cc$ into eq. (4.9), we can obtain
\begin{equation*}
  (s_\ast A+s_\ast A_se^{-i\omega_\ast s_\ast}-i\omega_\ast s_\ast\mathbf{I})E(T_2)=H'(\mathbf{I}+s_\ast A_se^{-i\omega_\ast s_\ast})\mathbf{c}-\delta H(A+A_se^{-i\omega_\ast s_\ast})\mathbf{c}
\end{equation*}
\begin{equation}\label{}
 \ \ \ \ \  +s_\ast MH^2\bar{H}.
\end{equation}
In order to get $E(T_2)$, there must exists $\mathbf{d}=(d_1,d_2,d_3)$ satisfying
 $$\mathbf{d}\times(s_\ast A+s_\ast A_se^{-i\omega_\ast s_\ast}-i\omega_\ast s_\ast\mathbf{I})=\mathbf{0}.$$
 \noindent
$\displaystyle$
 Then $d_1=\alpha d_2$ and $d_3=\beta d_2$ where
\begin{eqnarray*}
\begin{split}
  &\alpha=\frac{-\mathbf{v}_\ast}{r_1-2a_1r_1\mathbf{u}_\ast-b_1r_1\mathbf{v}_\ast e^{-i\omega_\ast s_\ast}-i\omega_\ast},
  \\
  &\beta=\frac{i\omega_\ast}{b_2r_2}.
  \end{split}
\end{eqnarray*}

 \noindent
$\displaystyle$
 To make $d_2$ unique, let

 $$\mathbf{d}(\mathbf{I}+s_\ast A_se^{-i\omega_\ast s_\ast})\mathbf{c}=1.$$

 \noindent
$\displaystyle$
 Then

 $$d_3=\frac{1}{\alpha c(1-b_1r_1s_\ast e^{-i\omega_\ast s_\ast}{v}_\ast)-\alpha b_1r_1{u}_\ast s_\ast e^{-i\omega_\ast s_\ast}+\beta+d}.$$
All the above conditions lead to the following equation
\begin{equation}
  H'=\delta H\Gamma_1-\Gamma_2H^2\bar{H},
\end{equation}
 \noindent
$\displaystyle$
where
 $\Gamma_1=\mathbf{d}(A+A_se^{-i\omega_\ast s_\ast})\mathbf{c}$ and $\Gamma_2=s_\ast\mathbf{d} M.$

  \noindent
$\displaystyle$
Equation (4.11) is the complex-value form of Hopf bifurcation at the critical value $s_\ast$.

We introduce the polar form $H=\rho e^{i\theta}$, it can be transformed into

 $$\rho'=\delta \rho \chi_1-\chi_2\rho^3,$$

 with $\chi_1=Re\Gamma_1$ and  $\chi_2=Re\Gamma_2$.

 According to \cite{Yuri}, it is easy to get the following theorem:

 \noindent
$\displaystyle$
 $\mathbf{Theorem 3.1}$  For system (1.2)

 (1)If $\chi_1\chi_2>0$, the Hopf bifurcation is supercritical and the bifurcating periodic solution is stable.

 (2)If $\chi_1\chi_2<0$, the Hopf bifurcation is subcritical and the bifurcating periodic solution is unstable.

Remark: If $\chi_1\chi_2>0$, there will be stable periodic oscillation  at the positive equilibrium point when it cross some critical value. It means the information system will have a oscillation at a short time and it will return back to the positive equilibrium again to keep a dynamic equilibrium. And if $\chi_1\chi_2<0$, the oscillation of information  may destroy network stability in online social networks and even cause a panic in the real society. Some appropriate measures should be taken by government to control the trend.
  \section{Numerical simulations}
  In this section, we will present a example to illustrate our theory results.

  Let $a_1=0.05$, $a_2=1.045$, $b_1=0.95$, $b_2=0.27$, $\mu=2$, $r=4$, $r_1=0.5$, $r_2=0.5$, then the system becomes
  \begin{equation*}
    \begin {cases}
    \begin{split}
    &\frac{{\rm d} {u}}{{\rm d} t}=0.5 {u}-0.025 {u}^2-0.475 {u}(t-s) {v}(t-s),\\
    &\frac{{\rm d} {v}}{{\rm d}t}=0.5 {v}-0.5225 {v}^2+0.135 {w},\\
    &\frac{{\rm d}{w}}{{\rm d}t}= {uv}-6 {w}.
    \end{split}
    \end{cases}
  \end{equation*}

  It is easy to know that there is only a positive equilibrium $(1,1,\frac{1}{6})$ and we can obtain $\omega_\ast=0.2004$, $s_\ast=2.015$.

  From theorem (3.1), we know that the system is locally asymptotically when $s=2$ [see Fig. 1(a)-(e)]. A Hopf bifurcation occurs when $s$ cross the critical value. When $s=2.02>s_\ast$, we can know that there are periodic oscillations from the positive equilibrium [see Fig. 2].

\begin{figure}[H]
\centering
\subfigure[]{
\label{Fig.sub.1}
\includegraphics[width=0.4\textwidth]{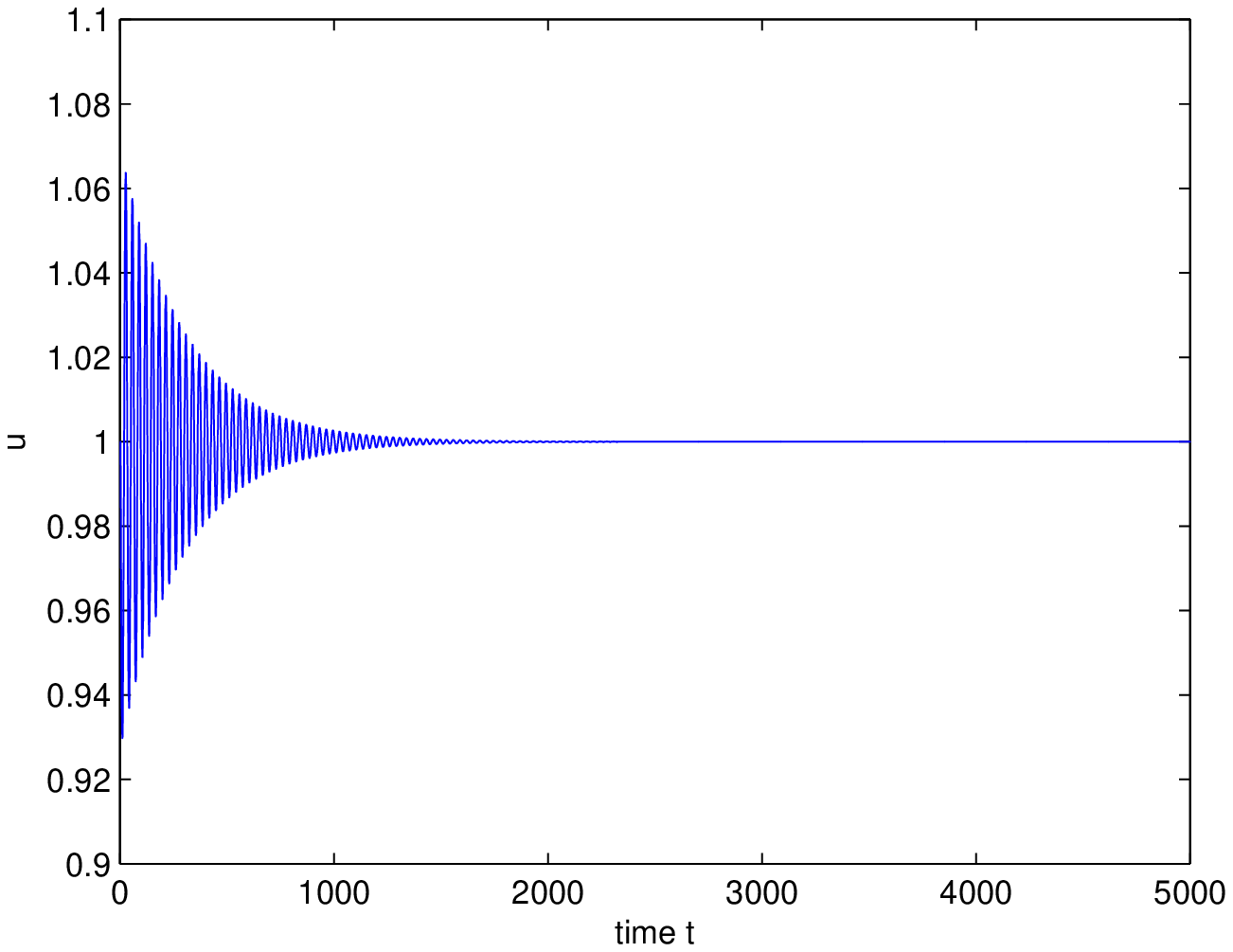}}
\subfigure[]{
\label{Fig.sub.2}
\includegraphics[width=0.4\textwidth]{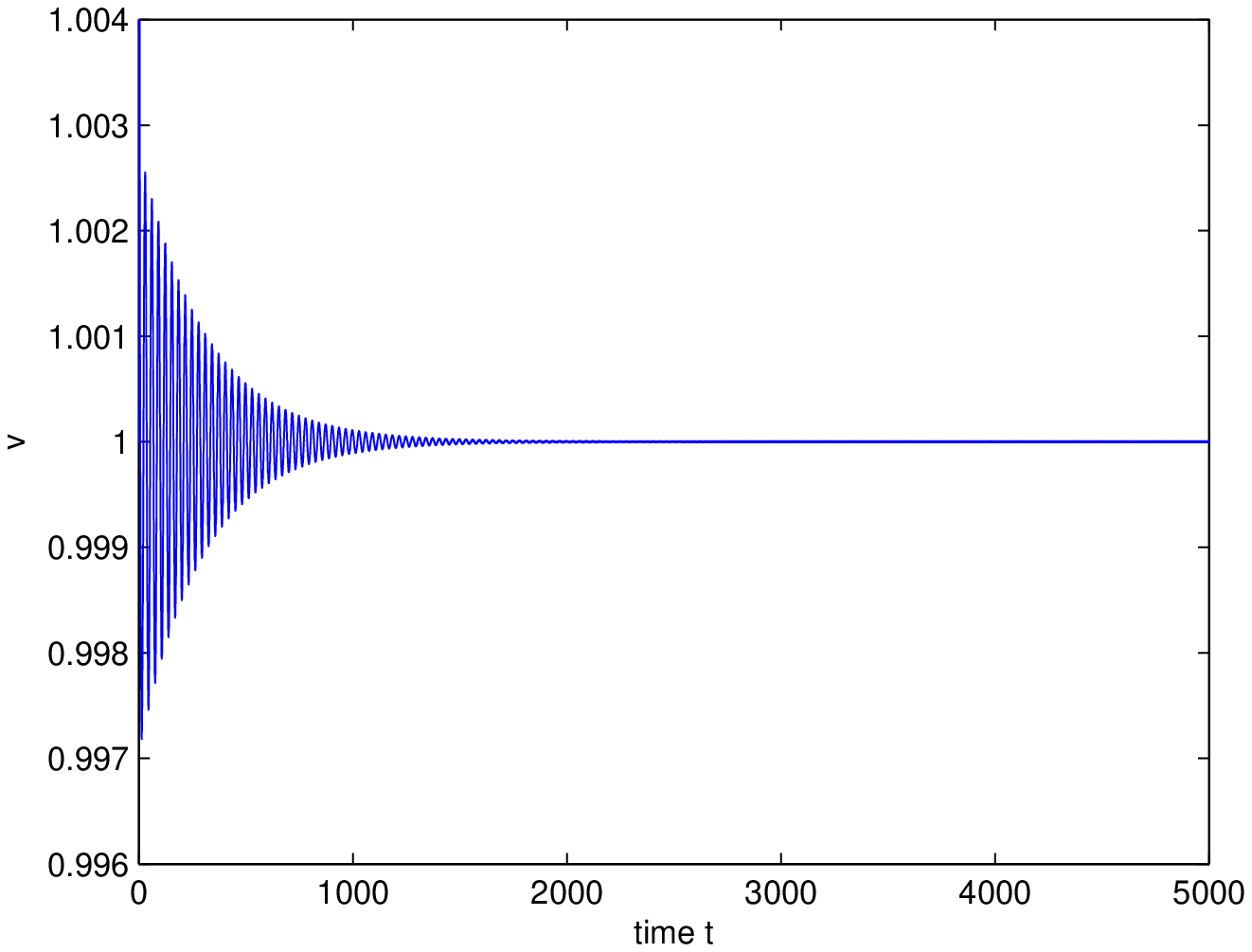}}
\label{Fig.lable}
\subfigure[]{
\label{Fig.sub.1}
\includegraphics[width=0.4\textwidth]{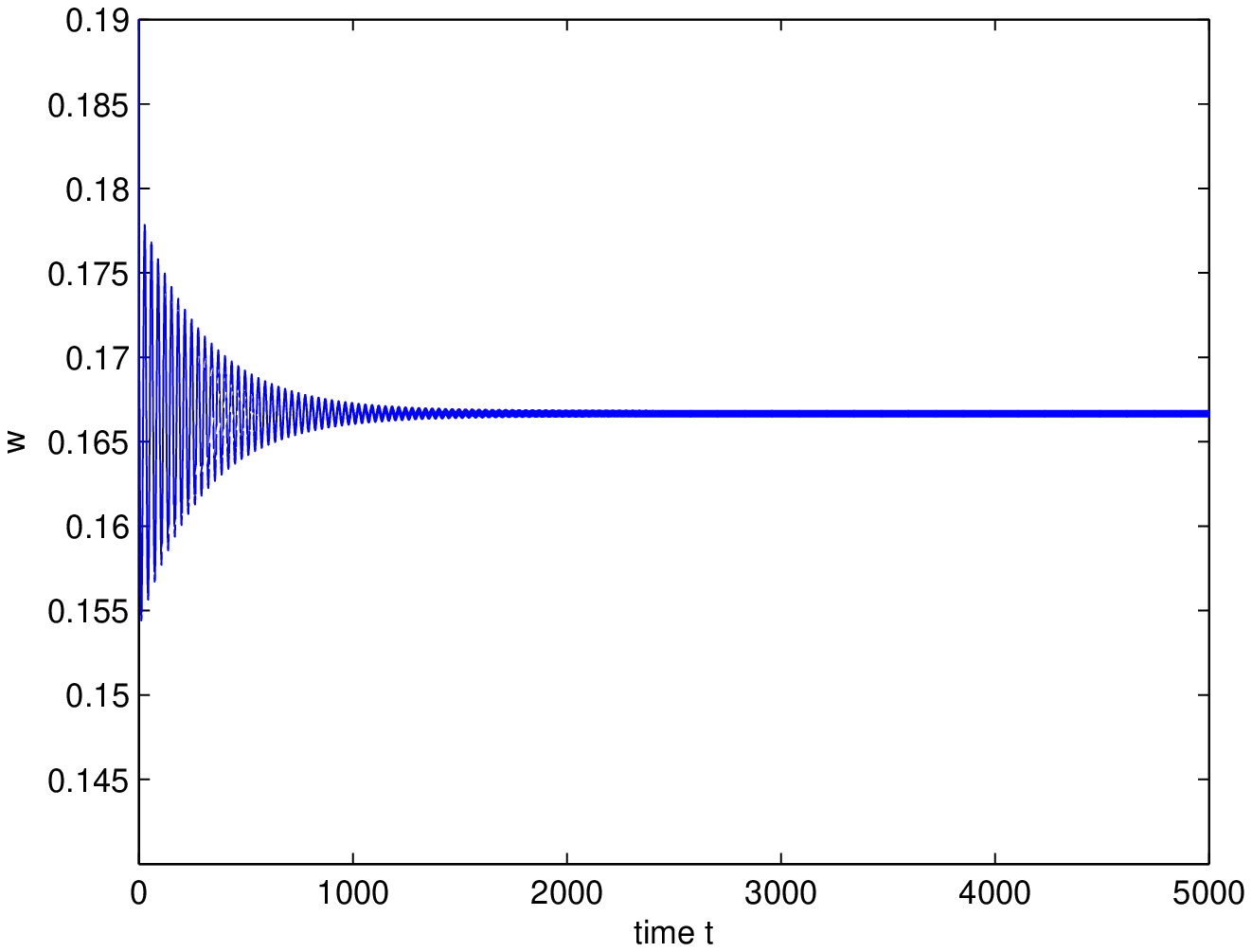}}
\subfigure[]{
\label{Fig.sub.1}
\includegraphics[width=0.4\textwidth]{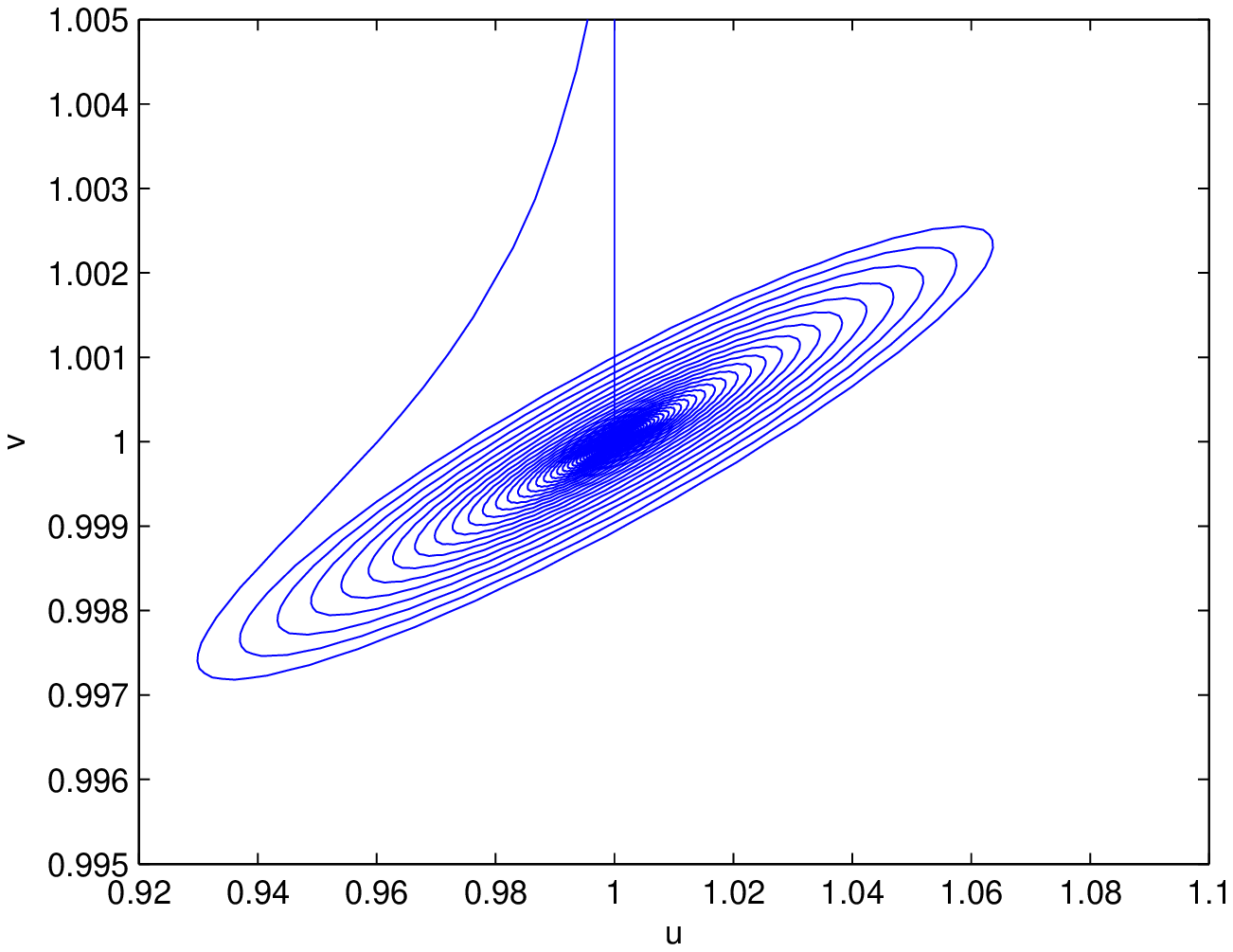}}
\label{Fig.lable}
\subfigure[]{
\label{Fig.sub.1}
\includegraphics[width=0.4\textwidth]{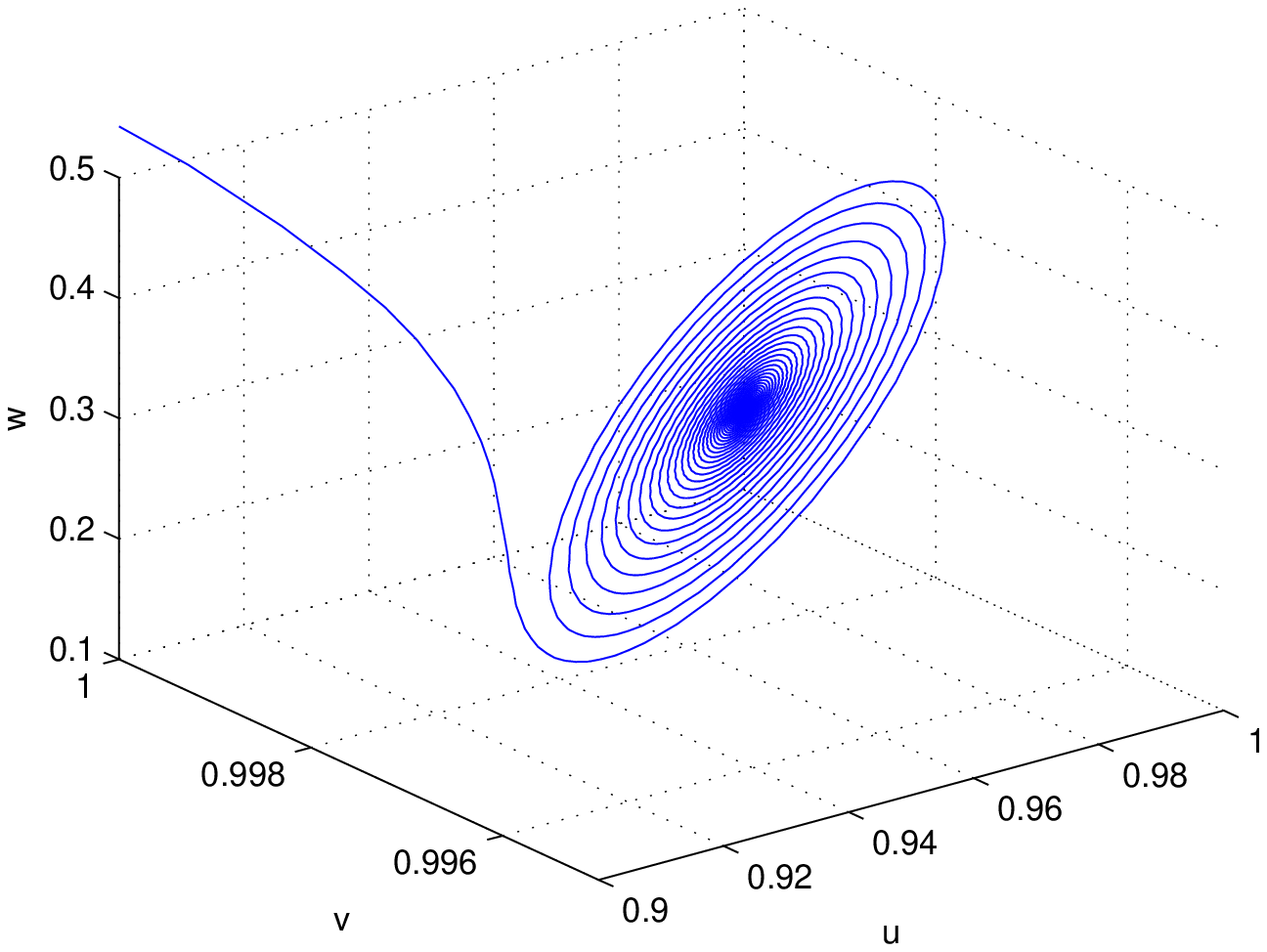}}
\caption{(a)-(c) show the waveform plots of the system and (d)-(e) show the phase portraits of the system. Let $a_1=0.05$, $a_2=1.045$, $b_1=0.95$, $b_2=0.27$, $\mu=2$, $r=4$, $r_1=0.5$, $r_2=0.5$ and $s=2<s_\ast=2.015$.}
\end{figure}

  \begin{figure}[H]
\centering
\subfigure[]{
\label{Fig.sub.1}
\includegraphics[width=0.4\textwidth]{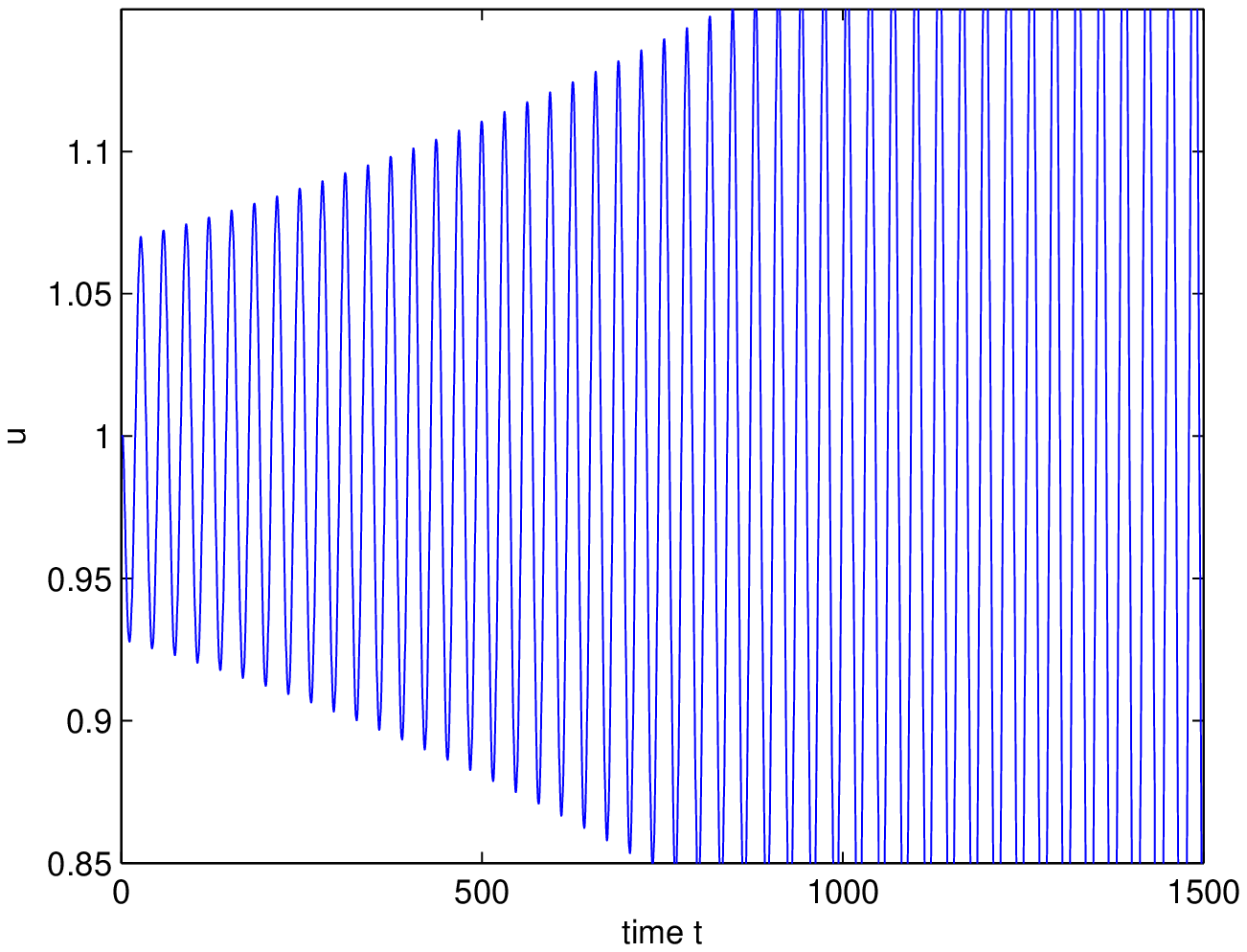}}
\subfigure[]{
\label{Fig.sub.2}
\includegraphics[width=0.4\textwidth]{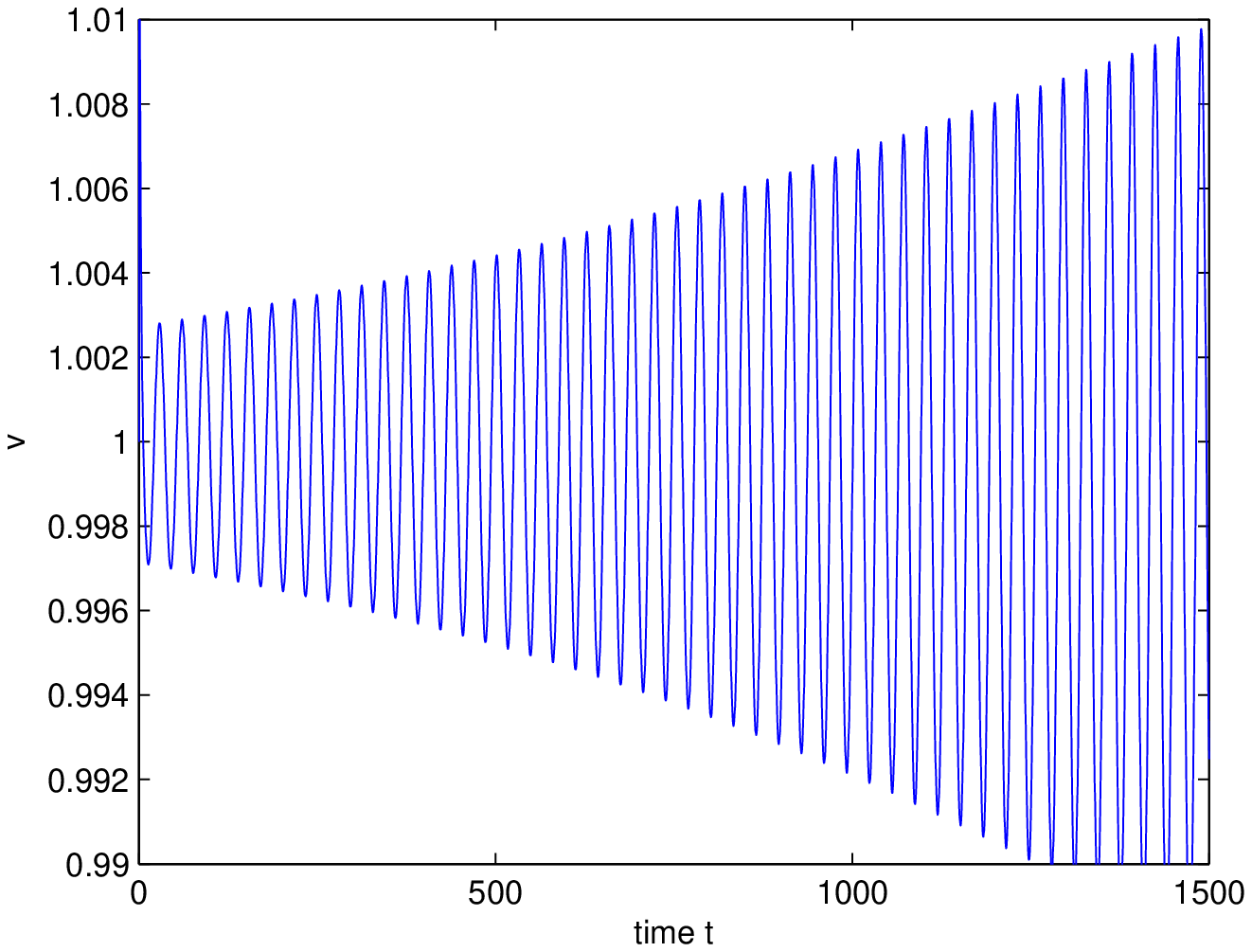}}
\label{Fig.lable}
\subfigure[]{
\label{Fig.sub.1}
\includegraphics[width=0.4\textwidth]{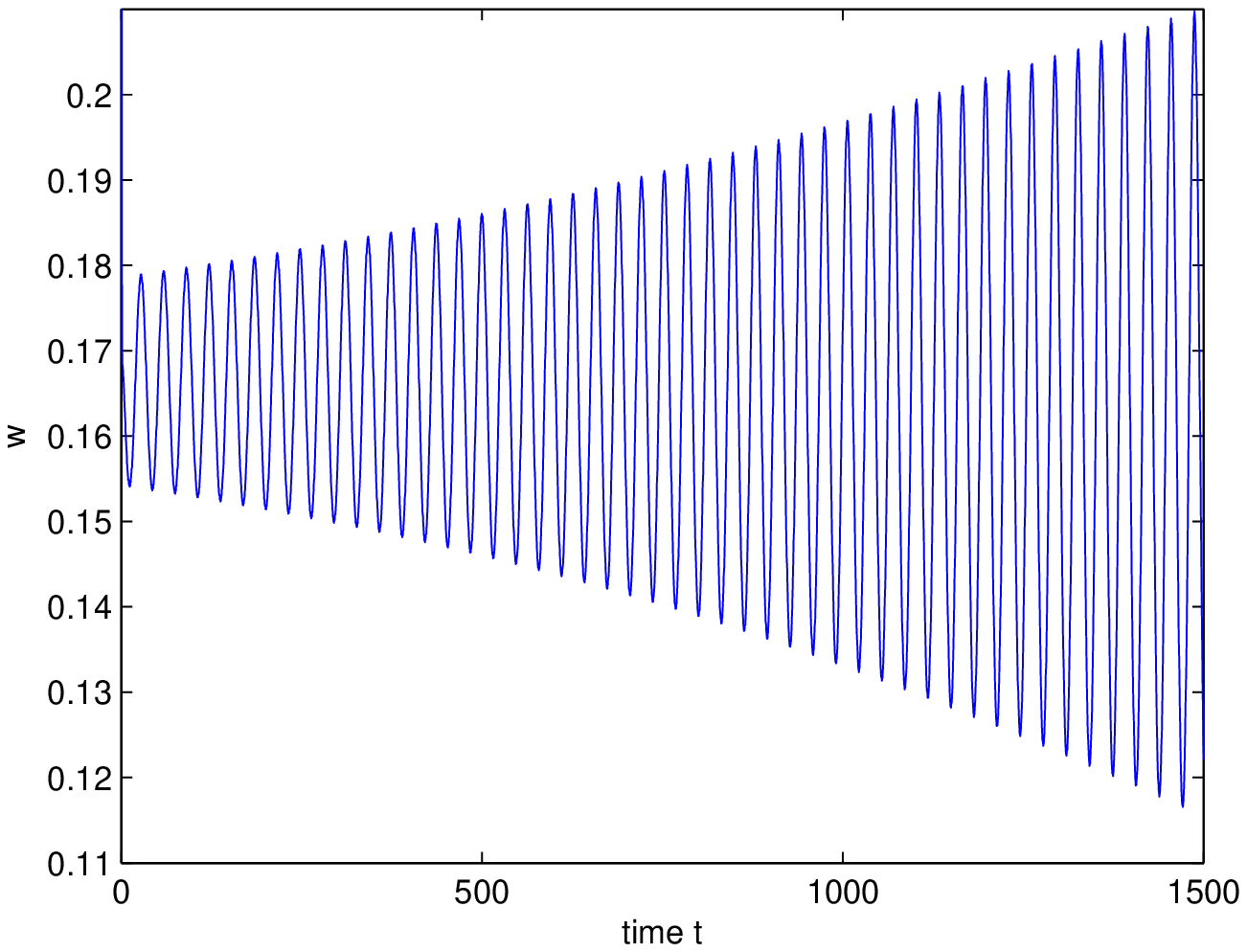}}
\subfigure[]{
\label{Fig.sub.1}
\includegraphics[width=0.4\textwidth]{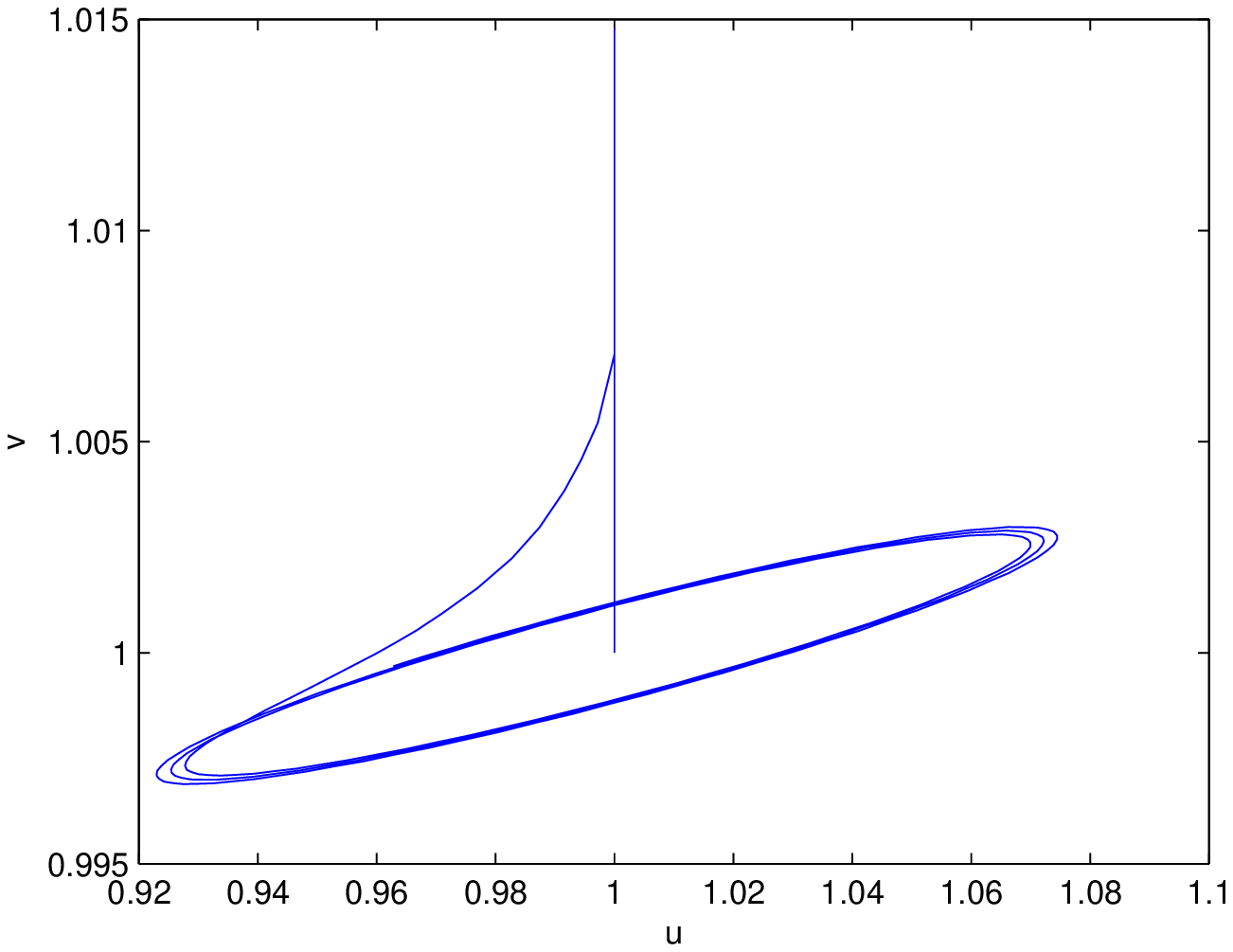}}
\label{Fig.lable}
\subfigure[]{
\label{Fig.sub.1}
\includegraphics[width=0.4\textwidth]{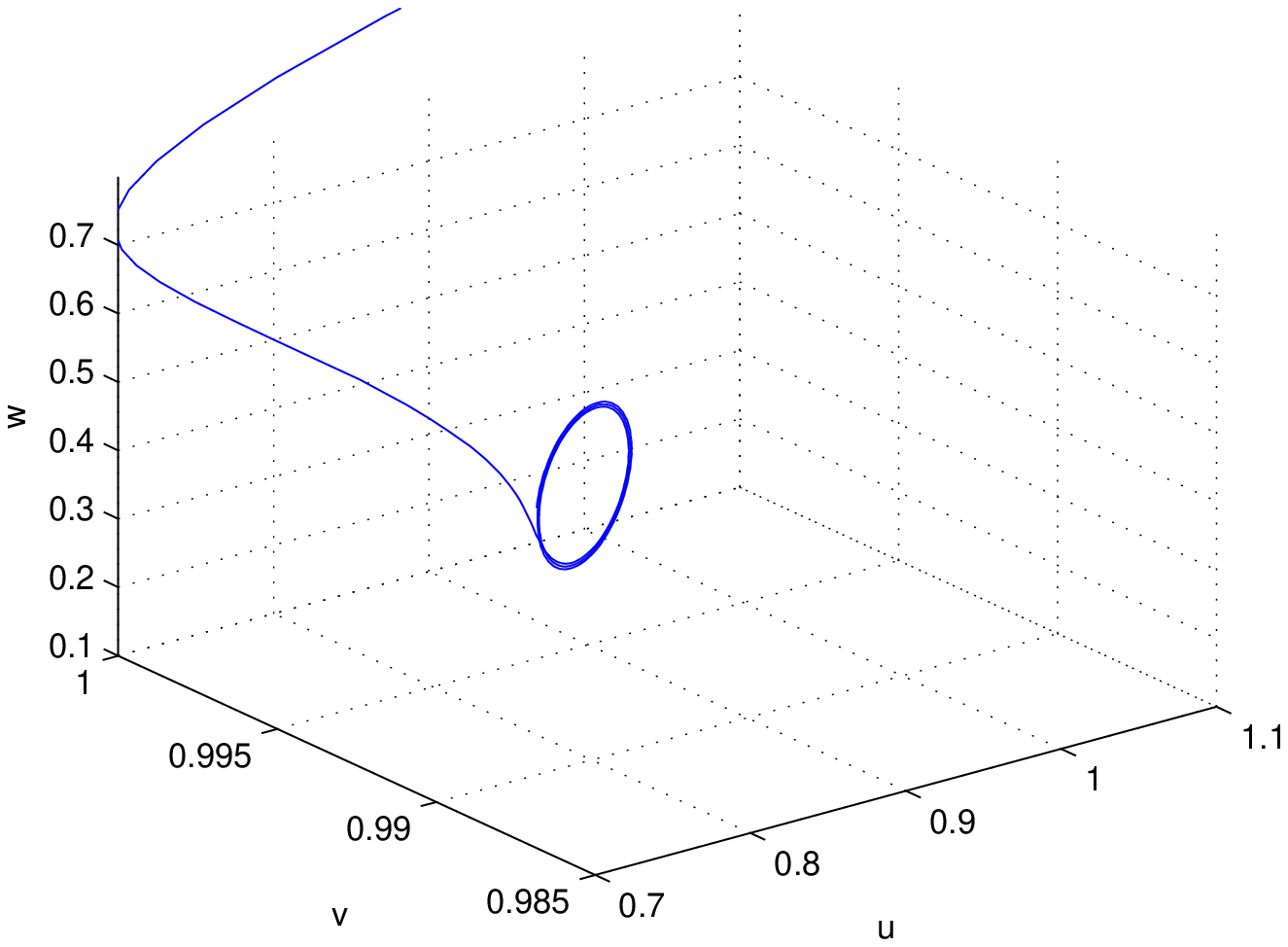}}
\caption{(a)-(c) show the waveform plots of the system and (d)-(e) show the phase portraits of the system. Let $a_1=0.05$, $a_2=1.045$, $b_1=0.95$, $b_2=0.27$, $\mu=2$, $r=4$, $r_1=0.5$, $r_2=0.5$ and $s=2.02>s_\ast=2.015$.}
\end{figure}

By computation, we can obtain that $\Gamma_1=2.9590-3.3311i$ and $\Gamma_2=(1.3813-5.0586i)\times10^{-3}$. It is obvious that $\chi_1\chi_2>0$, by the theorem (2.1), we know that the Hopf bifurcation is supercritical and the bifurcating periodic solution is stable, which can be seen in the numerical simulations.

\section{\Large{Acknowledgments}}
The research of the first author is supported by the the NSFC (Grant No. 11271339), the Plan for Scientific Innovation Talent of Henan Province (164200510011) and the ZDGD13001 Program.

  \end{document}